\magnification=\magstep1
\advance\voffset by 1 true cm
\baselineskip=15pt
\overfullrule =0pt 
\font\bbb=msbm10

\def\C{\hbox{\bbb C}}

\def\Z{\hbox{\bbb Z}}
\def\P{\hbox{\bbb P}}

\centerline {\bf ON THE MONODROMY OF COMPLEX POLYNOMIALS}

\vskip2truecm

\centerline {\bf by  Alexandru Dimca and Andr\'as N\'emethi}

\vskip2truecm

{\bf 1. Introduction}

\bigskip
Consider a polynomial function $f:\C^n \to \C$ with generic fiber $F$. Let $B_f$ be the bifurcation set of $f$, hence $f$ induces a smooth locally trivial fibration over $\C \setminus B_f$, see (2.1). Then, for any integer $q \geq 0$ and any coefficient ring $R$ there is an associated  monodromy representation
$$\rho (f) _q: \pi _1 (\C \setminus B_f,pt) \to Aut ({\tilde H}_q(F,R))$$
in (reduced) homology as well as a monodromy representation in (reduced) cohomology

$$\rho (f) ^q: \pi _1 (\C \setminus B_f,pt) \to Aut ({\tilde H}^q(F,R)).$$
Going around a large circle in $\C$ containing inside the bifurcation set gives rise to the so called monodromy operators at infinity, which we denote by $M_{\infty}(f)_q$ and  $M_{\infty}(f)^q$ respectively.

\bigskip

Our main result is the following.
\bigskip

\noindent {\bf Theorem A.}

{\it The following statements are equivalent.

\noindent (i) $M_{\infty}(f)_q=Id$;

\noindent (ii) the homology representation $\rho (f) _q$ is trivial.

If we assume that $R$ is a field or that $R =\Z$ and the integral (co)homology of $F$ is torsion free then (i) and (ii) are equivalent to the following:

\noindent (iii) $M_{\infty}(f)^q=Id$;

\noindent (iv) the cohomology representation $\rho (f) ^q$ is trivial.}

\bigskip

This result was first obtained in the case $n=2$ by the first author [D2], using  specific properties of this low dimensional situation plus some Hodge theory. The same case $n=2$ was also proved by Pi. Cassou-Nogu\`es and E. Artal-Bartolo [AC]  using the splice diagrams.

Here we prove it in full generality: actualy we obtained a much finer result, namely that the monodromy at infinity and a certain natural direct sum decomposition of the homology of $F$ in terms of vanishing cycles determine the monodromy representation, see Theorem (3.1). The role played by this decomposition is crucial, since there are examples of polynomials $\C^2 \to \C$ having distinct complex monodromy representations but whose monodromy operators at infinity have the same Jordan normal form, see the polynomials $f_B$ and $f_C$ in [D2].

The Theorem (3.1) is very unexpected since in many questions it is much easier to deal with the monodromy at infinity rather than with the whole of the monodromy representation. For instance, to compute the monodromy at infinity under the Thom-Sebastiani construction is not so difficult, see [N1], [N2], but to treat the corresponding monodromy representations is much harder, see [DN].

On the other hand, a lot of information on the monodromy at infinity is known,
at least in some special cases, see for instance [GN], and this can be hopefully used in conjunction with our Theorem (3.1) to get information on the monodromy representation itself.

As a by-product of Theorem A. we give conditions under which the generic fiber of the polynomial $f$ is acyclic in a certain degree $q$, generalizing in this way Corollary 2. from [D2]. Of course, in the case $n>2$, the Abhyankar-Moh-Suzuki Theorem is missing, so we cannot deduce from the acyclicity of the fibers that the polynomial $f$ is equivalent to a linear form, see Zaidenberg [Z] for such examples of polynomials.

\bigskip

\noindent {\bf Corollary B.} 

{\it Assume that $R=\Z$ and $q$ is a fixed nonnegative integer. Then the following
are equivalent:

\noindent (i) $M_{\infty}(f)_q=Id$ and ${\tilde H}_c^{2n-q-1}(F_b,R)=0$ for all $b \in B_f$;

\noindent (ii) ${\tilde H}_q(F,R)=0$ and $Ker(M_{\infty}(f)_{q-1}-Id)=0$.}

\bigskip

Here we use the notation $F_b=f^{-1}(b)$ for any $b \in \C$. For an $m$-dimensional affine variety $Y$, the reduced cohomology with compact supports $
{\tilde H}_c^{j}(Y,R)$ is defined to be just ${ H}_c^{j}(Y,R)$ when $j \not= 2m$ and ${\tilde H}_c^{2m}(Y,R)= Coker(p^*:  { H}_c^{2m}(\C^m,R) \to { H}_c^{2m}(Y,R)           $, where $p:Y \to \C^m$ is any finite morphism.

To simplify the notation we omit the coefficient ring $R$ when it plays no special role.

\vskip1truecm

{\bf 2. A Picard formula for some monodromy operators}

\bigskip

\noindent 
(2.1) Let $f: \C^n \to \C$ be a polynomial function. It is well known that there is a (minimal) finite bifurcation set $B_f$ in $\C$ such that $f$ is a $C^{\infty}$-locally trivial fibration over $\C \setminus B_f$.
If $c_0 \in \C$ is not in $B_f$, then $F=f^{-1}(c_0)$ is called the generic fiber of $f$.

The fundamental group $\pi _1(\C \setminus B_f, c_0)$ is a free group with $t=|B_f|$ generators. Note that in this paper we use the following convention: the product loop $a \cdot b$ is obtained by going first along the loop $b$ and then along the loop $a$.

Set $B_f=\{b_1,...,b_t\}$. Then a convenient way to construct a set of free generators $\gamma _i$ for this group is the following. First fix some paths $p_i:[0,1] \to \C$, $i=1,...,t$ such that we have

(i) $p_i(0)=c_0$, $p_i(1)=b_i$;

(ii) (image $p_i) \cap$ (image $p_j)=\{c_0\}$ for any $i \not= j$;

(iii) the paths $p_1$,...,$p_t$ leave $c_0$ in the counterclockwise order.
\bigskip

\noindent (2.2) {\bf Definition.} \ 
A family of paths $(p_i)_{i=1,...,t}$ with the properties (i), (ii) and (iii)
is called a star with base point $c_0$ and end points $b_1$,...,$b_t$. 

\bigskip

Once we fix a star as above, we have automatically a base $(\gamma _i)_{i=1,...,t}$ of $ \pi _1(\C \setminus B_f, c_0)$. Indeed, consider small discs $D_i$ centered at $b_i$ such that (image $p_i) \cap D_i$ is a topological radius in $D_i$. Set (image $p_i) \cap \partial D_i=c_i=p_i(t_i)$. Then ${\tilde p}_i:[0,t_i] \to \C$, ${\tilde p}_i(t)=p_i(t)$ connects $c_0$ to $c_i$. Let $\partial D_i$ be the loop with both end points $c_i$ and going once around $b_i$ counterclockwise. 

Then $\gamma_i=({\tilde p}_i)^{-1}\circ  \partial D_i \circ {\tilde p}_i$ for $i=1,...,t$ are loops based at $c_0$ and their homotopy classes (still denoted by $\gamma_i$) generate the group $\pi _1(\C \setminus B_f, c_0)$. We emphasize that the ordered set of homotopy classes $\gamma_1,...,\gamma_t$ is canonically associated with the star $p_1,....,p_t$ and it is independent on any other choices.

\bigskip

\noindent 
(2.3) {\it The sum decomposition of ${\tilde H}_q(F)$ associated with a star.}

\bigskip

Consider as above a polynomial $f:\C^n \to \C$ with a bifurcation set $B_f$,
a point $c_0$ not in $B_f$ and a star with base point $c_0$ and end points in $B_f$. This provides a direct sum decomposition of ${\tilde H}_q(F)$, see [B], [DN], [NN1], [ST1], [ST2] for various degrees of generality.

Indeed, first we have the obvious isomorphism 
$\partial ^{-1}: {\tilde H}_{q}(F) \to H_{q+1}(\C^n,F)$, next we can replace the second group by $H_{q+1}(f^{-1}(\cup_i p_i([0,1]) \cup \cup_i D_i), F)$ using a deformation retract and finally by excision we may replace the last group by $\oplus_i H_{q+1}(f^{-1}(D_i),F_{c_i})$. Next we can define a {\it subgroup of vanishing cycles} $V_{q,i}$ in ${\tilde H}_{q}(F)$ for $i=1,...,t$ by pushing the subgroup $\partial H_{q+1}(f^{-1}(D_i),F_{c_i})$
in ${\tilde H}_{q}(F_{c_i})$ to ${\tilde H}_{q}(F)$ along the path $p_i$.

This sum decomposition
$$\leqno {(*_q)} ~~~~~~~~~~~~~~~~~~~~~ {\tilde H}_{q}(F)= \oplus_{i=1}^t\,V_{q,i}$$
depends essentially on the choice of the star.

If $R \to K$ is a ring extension with $K=R$ or $K$ a field, then by tensor product by $K$ over $R$ in $(*_q)$ we get the corresponding direct sum decomposition of ${\tilde H}_{q}(F,K)$
which will be denoted sometimes by $(*_q(K))$.

Some notation: let 
$m_j=\rho (f)_q (\gamma_j)={\tilde H}_q(T_j)   \in Aut ({\tilde H}_q(F))$ and  
$m_j^*=\rho (f)^q (\gamma_j)= {\tilde H}^q(T_j)^{-1} \in Aut ({\tilde H}^q(F))$ 
be the monodromy operators in homology and in cohomology
associated to a star as above. Here $T_j: F \to F$ is the homeomorphism obtained by parallel transport along the path $\gamma_j$.

It is clear that these operators completely determine the corresponding monodromy representation and that $m_k \cdot ... \cdot m_1 =M_{\infty}(f)_q$ and
$m_k^* \cdot ... \cdot m_1^* =M_{\infty}(f)^q$

When we need to keep track of the degree in which these operators act, we write
$m_{q,j}$ instead of $m_j$ and $m_j^q$ instead of $m_j^*$.

\bigskip

The following result provides key information on these monodromy operators in homology (exactly as in the case of isolated hypersurface singularities and/or the classical  theory of Lefschetz pencils) and can be found in the global setting in various recent preprints: [DN], Lemma (2.4), [NN1], Theorem (1.4), [NN2], Theorem (3.1) and [ST2], Proposition (2.2). The proof of it being very simple, we reproduce it below.

\bigskip

\noindent (2.4) {\bf Lemma.} \ {\it 
For any cycle $[c] \in {\tilde H}_{q}(F) $:}
$$(m_j([c])-[c]) \in V_{q,j}.$$

\bigskip

\noindent {\it Proof.} \
The geometric monodromy along the path $\gamma _j$ is defined as follows. 
Take a $q$-cycle $c$ and move it along the path $\gamma _j$ to get a cycle 
$m_{geom}(c)$ in $F$. This means that there is a $(q+1)$-cycle $d$ such that 
$ \partial d=m_{geom}(c)-c$ and $f(d) \subset \gamma_j([0,1])$. But then 
$m_j([c])-[c]=\partial [d]\in H_{q,j}$. 
%Note that we also have $(m_j^{-1}([c])-[c]) \in H_{*,j}$
%with the same proof.

\bigskip

\noindent (2.5)\ The above Lemma shows that in the sum-decomposition 
$\oplus_{i=1}^t\, V_{q,i}$, the monodromy $m_j=\rho_f(\gamma_j)$ 
($1\leq j\leq t$) has a 
block-decomposition  of the following form:

$$ \leqno {(**_q) } ~~~~~~~~~~~~~~~~~~~~  m_j= 
\ \pmatrix{1& &&&& \cr
\ & \ddots &&&& \cr
m_{j1}&\cdots& m_{jj}&m_{j,j+1}& \cdots& m_{jt}\cr
\ &&&1&& \cr
\ &&&& \ddots& \cr
\ &&&&& 1\cr},$$
where all the entries, excepting the diagonal and the $j$-th line, are trivial.
Here $m_{ji}:V_{q,i}\to V_{q,j}$ is given by $m_{ji}(a_i)=m_j(a_i)-a_i$ (for
$i\not=j$); and $m_{jj}:V_{q,j}\to V_{q,j}$ is $m_{jj}
(a_j)=m_j(a_j)$. 

\bigskip 

\noindent (2.6) {\bf Remarks.} \ 
(i) Lemma (2.4) is a weak global affine Picard-Lefschetz formula 
(corresponding to the result of Picard) 
and it says that $m_j([c])-[c]$ can be localized in $f^{-1}(D_j)$.

(ii) We have  similar decompositions $(*_q)$ and $(**_q)$ in the more general case of a regular function $f:X \to \C$ where $X$ is an algebraic varity with $H_{q+1}(X)={\tilde H}_{q}(X)=0$. Hence all our results below hold in this more general setting.

(iii) Note that in some cases the decomposition (2.3) is completely determined by the monodromy representation, see Lemma (8.1) in [DN].

\vskip1truecm

{\bf 3. The monodromy at infinity determines the monodromy representation}

\bigskip

In this section we prove the following result which implies Theorem A.
\bigskip

\noindent {\bf Theorem (3.1)}

{\it Let $R=\Z$ and fix a non-negative integer $q$ and an extension $\Z \to K$. Then:

\noindent (i) The monodromy operator at infinity $M_{\infty}(f)_q$ and the decomposition $(*_q)$ completely determines the monodromy representation $\rho (f)_q$.

\noindent (ii) If $v \in  {\tilde H}_{q}(F,K) $ has components $(v_1,...,v_t)$ with respect to the decomposition $(*_q(K))$ then
$$M_{\infty}(f)_q(v)=av$$
for some $a \in K$ if and only if
$$m_k \cdot ...\cdot m_1 (v)=(av_1,...,av_k,v_{k+1},...,v_t)$$
for all $k=1,2,...,t$.
In particular
$$Ker(M_{\infty}(f)_q-1)=  {\tilde H}_{q}(F)^{inv} $$
the invariant homology under the monodromy representation. }

\bigskip

The precise meaning of $(i)$ above is that if we know the operator $M_{\infty}(f)_q$ as a matrix of blocks of type $(**_q)$, then this gives the matrices corresponding to all the operators $m_k$, $k=1,...,t.$

Before giving the proof we would like to emphasize some points.

\bigskip

(a) The first claim $(i)$ holds for the cohomology monodromy representation as well  when $H_*(F,R)$ is torsion free. Then  $M_{\infty}(f)_q$ is the inverse transpose of $M_{\infty}(f)^q$ in the following sense. The decomposition $(*_q)$ induces a similar decomposition for $
{\tilde H}^{q}(F)$, see [DN], (8.2), which we denote by $(*^q)$. If the matrix $A$ of $M_{\infty}(f)^q$ is known with respect to this decomposition $(*^q)$, then the matrix of $M_{\infty}(f)_q$ with respect to the decomposition $(*_q)$ is just the inverse transpose matrix $(A^t)^{-1}$.

\bigskip

(b) The second part of (3.1) is definitely false for cohomology, e.g. for $R$ a field, there are examples when
$$dim Ker(M_{\infty}(f)^q-1) >  dim{\tilde H}^{q}(F)^{inv} $$
see for instance the Brian\c con polynomial in [D2].

More precisely we have the following situation when $R$ is a field. For a single monodromy operator we have an obvious equality
$$dim Ker(m_k-1)=dim Ker(m_k^*-1)$$
and various formulas for this dimension have been obtained in [ACD] and, in full generality, in [NN1]. Moreover the vector spaces $( Ker(m_k^*-1))_k$ are in general position in ${\tilde H}^{q}(F)$, see [ACD] and [NN1]. This implies the following inequality (which is strict in general):
$$dim{\tilde H}_{q}(F)^{inv}=dimKer(M_{\infty}(f)_q-1)=dimKer(M_{\infty}(f)^q-1) \geq dim{\tilde H}^{q}(F)^{inv}.$$
This discussion shows that homology and cohomology monodromy representations have each their advantages and hence it is usefull to study and use them both.

\bigskip

\noindent {\it (3.2) PROOF of Theorem (3.1)}

\bigskip

The shape of the matrices in $(**_q)$ implies the following obvious and equivalent facts. For any $k=1,...,t$ we have

$(\cdot)$ The operator $m_k$ is determined by the composition $p_k \cdot m_k$ where $p_k: 
{\tilde H}_{q}(F) \to V_{q,k}$ is the projection onto the $k$-th factor in $(*_q)$.

$(\cdot \cdot)$ $m_k(v_1,...,v_t)=(v_1,...,v_{k-1},v'_k,v_{k+1},...,v_t)$, i.e. the operator $m_k$ acts only on the $k$-th component of any vector $v$.

The formula
$$\leqno {(***)} ~~~~~~~~~~~~~~~~~~~~ m'' \cdot m_k \cdot m'= M_{\infty}(f)_q$$
 where $m'=m_{k-1} \cdot ...\cdot m_1$ and  $m''=m_t \cdot ... \cdot m_{k+1}$
implies that

$$  p_k \cdot m_k =p_k \cdot M_{\infty}(f)_q \cdot m'^{-1}.$$
Indeed, the operator $m''$ leaves the first $k$ components of a vector $v$ unchanged. But this shows that, if $m_1$,...,$m_{k-1}$ are already determined by 
$M_{\infty}(f)_q$, then the same holds for $m_k$, i.e. (3.1) $(i)$ follows by induction on $k$.

To get the second part $(ii)$, we use again the decomposition $(***)$ and the property of $m''$, to get that $M_{\infty}(f)_q(v)=av$ implies
$$m_k \cdot m' (v)=(av_1,...,av_k,v_{k+1},...,v_t).$$
Note for this that $m_k \cdot m'$ leaves unchanged the last $(t-k)$ components of any vector.

\bigskip

\noindent {\it (3.3)  Proof of Corollary B.}

\bigskip

This proof uses our Therem A. and the following exact sequence from [NN1], Theorem 1.4.

$$\leqno {(E)}: ~~~~~~~~~~~~~~0 \to Im(m_{q,j}-1) \to V_{q,j} \to {\tilde H}_c^{2n-q-1}(F_{b_j}) \to Ker(m_{q-1,j}-1) \to 0. $$
In fact, in {\it loc.cit} the group ${\tilde H}_c^{2n-q-1}(F_{b_j})$ is replaced by the group $H^{2n-q-1}(Y,\partial Y)$ where $Y$ (resp. $\partial Y$) is the intersection of the special fiber $F_{b_j}$ with a very large closed ball $ B$ in $\C^n$ (resp. with the boundary of this ball).
The two groups are clearly isomorphic and we prefer to work with the cohomology with compact supports.

If $S_b$ is the singular part of the fiber $F_b$ and $d=dim(S_b)$, then the natural morphism
$$\leqno{(3.4)} ~~~~~~~ H_c^j(F_b \setminus S_b) \to H_c^j(F_b)$$
is surjective for $j=2d+1$ and an isomorphism for $j>2d+1$.
In particular for $R$ a field:

$(i)$ $dimH_c^{2n-2}(F_b)$ is equal to the number of irreducible components of $F_b$ as in [NN1];

$(ii)$ $dimH_c^{2n-3}(F_b)= dim H^1(F_b \setminus S_b)$ if $d<n-2.$
\bigskip

To end the proof, note that  the condition (i) (resp. (ii)) in Corollary B. is equivalent to the vanishing of the first and the third (resp. the second and the fourth) terms in the exact sequence $(E)$.

\vskip1truecm

{\bf 4. Some examples}

\bigskip

 \noindent {\bf (4.1) Example} (non-trivial monodromy on torsion in homology)

Here we give an example of a polynomial whose generic fiber $F$ has torsion in integral homology and the action of the monodromy on this torsion is non-trivial. Hence here we take $R=\Z$.

Let $Q$ be a three-cuspidal quartic  curve in the projective plane $\P^2$. In an appropiate coordinate system on $\P^2$, the curve $Q$ is the zero set of the polynomial
$$ f(x,y,z)=x^2y^2+y^2z^2+x^2z^2-2xyz(x+y+z).$$
It is well known that the fundamental group $\pi _1(\P^2 \setminus Q)$ is a finite non-abelian group of order 12, see for instance [D1], p.131.

Regarding now $f$ as a polynomial function $\C^3 \to \C$ we have

(i) $t=1,~~~b_1=0$;

(ii) the only monodromy homeomorphism $h:F \to F$ can be represented by
$$h(x,y,z)=(ix,iy,iz)$$
where $i^2=-1$, hence $h^4=1$.

(iii) If $<h>$ denotes the cyclic group of order 4 spanned by $h$, then the quotient space $F/<h>$ is just $\P^2 \setminus Q$. This implies that $H_1(F)=\pi _1(F)=\Z/3\Z$, see also Corollary (4.14) in [D1], p.133.
\bigskip

In the exact sequence $(E)$ we have in this case:
$V_{1,1}=H_1(F)=\Z/3\Z$,  ${\tilde H}_c^{4}(F_0)=0$ since $F_0$ is irreducible
and hence $Im(m_{1,1}-1)=\Z/3\Z$. In other words, the monodromy operator in homology satisfies
$$m_{1,1}=M_{\infty}(f)_1 \not= 1$$
On the other hand $H^1(F)=0$ and so $M_{\infty}(f)^1=1$. This example explains our restrictions on torsion in Theorem A.

\bigskip

The fact that $H_1(F)$ is torsion was already remarked by D. Siersma in [S], Example 7.4. The methods developed in his paper can also be used to show that 
in this case $m_{1,1}=M_{\infty}(f)_1 \not= 1$.
\bigskip

\noindent {\bf (4.2) Example} (trivial monodromy versus trivial intersection form)

It was shown by Neumann and Norbury [NN2] that the total space of the fibration
$ f: f^{-1}(S^1_r) \to S^1_r$ for $r>>0$ which gives the monodromy at infinity for the polynomial $f$ can be embedded in a natural way as an open subset of a large sphere $S^{2n-1}$ in $\C^n$.

This gives a Seifert form $L:{\tilde H}_{n-1}(F) \times {\tilde H}_{n-1}(F) \to \Z$
defined in the usual way as well as the usual relation
$$ L(1-M_{\infty}(f)_{n-1})=S$$
where $S$ is the intersection form on ${\tilde H} _{n-1}(F)$.
This equality gives the following.
\bigskip

\noindent {\it (4.3) $M_{\infty}(f)_{n-1}=1$ implies $S=0$, i.e.
$S(x,y)=0$ for any $x,y \in {\tilde H}_{n-1}(F)$.}

\bigskip

The converse implication is false in general, since the Seifert form $L$ can be degenerate. For instance the polynomial $f=x^2y^2+y: \C^2 \to \C$ has $M_{\infty}(f)_1 \not=1$ and
$S=0$ since $F$ is homeomorphic to $\C \setminus \{0,1\}$, see Bailly-Maitre [BM] for details.

This converse implication is nevertheless true for 
 the M-tame polynomials introduced by the second author and A. Zaharia in [NZ1], as a generalization of the tame polynomials introduced by Broughton in [B].

In fact, it was shown in [NZ2] that for an M-tame polynomial the fibration
$ f: f^{-1}(S^1_r) \to S^1_r$ for $r>>0$ which gives the monodromy at infinity for the polynomial $f$ is equivalent to the Milnor fibration at infinity
$\phi : S^{2n-1}_r \setminus f^{-1}(0) \to S^1$, $\phi (x)=f(x)/|f(x)|$.
This implies in the standard way that in this case the Seifert form is non-degenerated.

\bigskip

There is an analogy between our Example (4.2) and the  relations between $Ker(M(f)_*-1)$ and $KerS$ in the case of a local non-isolated hypersurface singularity as discussed in detail by D. Siersma in [S], section 7. Here $M(f)_*$ stands for the corresponding local monodromy operators.

\vskip1truecm

\noindent {\bf REFERENCES}
\bigskip

\item{[ACD]} E. Artal Bartolo, P. Cassou-Nogu\`es and A. Dimca: Sur la topologie des polyn\^omes complexes, Progress in Math. 162, Birkh\"auser 1998, pp 317-343.

\item{[AC]}  E. Artal Bartolo and P. Cassou-Nogu\`es:    
 Polyn\^ome d'Alexander \`a l'infini d'un polyn\^ome \`a deux variables,  Revista Matem\'atica de la Universidad Complutense de Madrid (to appear).

\item{[BM]} G. Bailly-Maitre: Sur le syst\`eme local de Gauss-Manin d'un polyn\^ome de deux variables, Bull. Soc. Math. France (to appear).

\item{[B]} S. A. Broughton: Milnor numbers and the topology of polynomial hypersurfaces, Invent. Math. 92(1988), 217-241.

\item{[D1]} A. Dimca: Singularities and Topology of Hypersurfaces,
Universitext, Springer, 1992.

\item{[D2]} A. Dimca: Monodromy at infinity for polynomials in two variables, J. Algebraic Geometry 7(1998) 771-779.

\item{[DN]} A. Dimca and A. N\'emethi: Thom Sebastiani construction and monodromy of polynomials, preprint 98 (March 1999), Bordeaux University.

\item{[GN]} R. Garc\'\i a L\'opez and A. N\'emethi: 
On the monodromy at infinity of a polynomial map, I. Compositio Math.,
100 (1996), 205-231; II.  Compositio Math. 115(1999), 1-20.

\item {[N1]} A. N\'emethi: Global Sebastiani-Thom theorem for polynomial maps, J. Math. Soc. Japan 43(1991), 213-218.

\item {[N2]} A. N\'emethi: Generalized local and global Sebastiani-Thom type theorems, Compositio Math. 80(1991), 1-14.

\item {[N3]} A. N\'emethi: On the Seifert form at infinity associated with polynomial maps, J. Math. Soc. Japan 51(1999),63-70.

\item{[NZ1]} A. N\'emethi and A. Zaharia: On the bifurcation set of a
polynomial function and Newton boundary, Publ. RIMS Kyoto
Univ. 26(1990), 681-689.

\item{[NZ2]} A. N\'emethi and A. Zaharia: Milnor fibration at infinity, Ingag. Math. 3(1992), 323-335.

\item{[NN1]} W. Neumann and P. Norbury: Monodromy and vanishing cycles of complex polynomials, Duke Math. J., to appear.

\item{[NN2]} W. Neumann and P. Norbury: Unfolding polynomial maps at infinity, preprint (October 1999), Melbourne University. math.AG/9910054.

\item{[S]} D. Siersma: Variation mappings on singularities with a 1-dimensional
critical locus, Topology 30(1991), 445-469.

\item{[ST1]} D. Siersma, M. Tib\u ar: Singularities at infinity and their vanishing cycles, Duke Math. J. 80(1995), 771-783.

\item{[ST2]} D. Siersma, M. Tib\u ar: Vanishing cycles and singularities of meromorphic functions, preprint no. 1105 (May 1999), Utrecht University. math.AG/9905108.

\item{[Z]} M. Zaidenberg: Lectures on exotic algebraic structures on affine spaces, Schriftreihe des Graduierten Kollegs Geometrie und Mathematische Physik, Heft 24, Ruhr-Universit\" at Bochum (1997).

\bigskip

Laboratoire de Math\'ematiques Pures de Bordeaux

Universit\'e Bordeaux I

33405 Talence Cedex, FRANCE

email: dimca@math.u-bordeaux.fr

\bigskip

Department of Mathematics,

Ohio State  University,

Columbus, Ohio 43210, USA

email: nemethi@math.ohio-state.edu

\bye